\documentclass[11pt,twoside]{amsart}

\theoremstyle{plain}
\newtheorem{thm}{Theorem}[section]

\author{Valery Alexeev} 
\address{Department of Mathematics\\
University of Georgia\\
Athens, GA 30602}
\email{valery@math.uga.edu}

\newcommand{\bX}{{\mathbb X}}

\newcommand{\bT}{{\mathbb T}}
\newcommand{\bP}{{\mathbb P}}
\newcommand{\bQ}{{\mathbb Q}}
\newcommand{\bZ}{{\mathbb Z}}

\newcommand{\bR}{{\mathbb R}}
\newcommand{\bC}{{\mathbb C}}
\newcommand{\bG}{{\mathbb G}}

\newcommand{\cR}{{\mathcal R}}
\newcommand{\cL}{{\mathcal L}}

\newcommand{\cM}{{\mathcal M}}

\newcommand{\cP}{{\mathcal P}}

\newcommand{\cO}{{\mathcal O}}

\newcommand{\cX}{{\mathcal X}}

\newcommand{\cS}{{\mathcal S}}

\newcommand{\fm}{\mathfrak m}

\newcommand{\hc}{\widehat c}

\newcommand{\uX}{\underline{X}}
\newcommand{\uY}{\underline{Y}}

\newcommand{\wG}{\widetilde G}

\newcommand{\wV}{\widetilde V}
\newcommand{\wP}{\widetilde P}

\newcommand{\wdL}{\widetilde {\mathbf L}}

\newcommand{\skipline}{\mbox{}\smallskip}

\DeclareMathSymbol{\curvearrowright}{\mathrel}{AMSb}{"79}
\DeclareMathSymbol\rightsquigarrow {\mathrel}{AMSa}{"20}
\DeclareMathSymbol\square {\mathord}{AMSa}{"03}
\DeclareMathSymbol{\ltimes}         {\mathbin}{AMSb}{"6E}
\DeclareMathSymbol{\nmid}           {\mathrel}{AMSb}{"2D}
\DeclareMathSymbol{\twoheadrightarrow}  {\mathrel}{AMSa}{"10}

\newcommand{\isoto}{{\overset{\sim}{\rightarrow}}}
\newcommand{\ratmap}{- \kern -3pt \to}

\newcommand{\into}{\hookrightarrow}

\newcommand{\follows}{\Rightarrow}

\newcommand{\emb}{\operatorname{emb}}

\newcommand{\orb}{\operatorname{orb}}

\newcommand{\Cone}{\operatorname{Cone}}

\newcommand{\id}{\operatorname{id}}

\newcommand{\GL}{\operatorname{GL}}

\newcommand{\Proj}{\operatorname{Proj}}

\newcommand{\Spec}{\operatorname{Spec}}

\theoremstyle{plain}
\newtheorem{question}[thm]{Question}

\newtheorem{lem}[thm]{Lemma}
\newtheorem{cor}[thm]{Corollary}

\theoremstyle{definition}

\newtheorem{defn}[thm]{Definition}
\newtheorem{setup}[thm]{Setup}

\newtheorem{notation}[thm]{Notation}

\newtheorem{assume}[thm]{Assumption}

\newtheorem{saynum}[thm]{}
\newtheorem{exmp}[thm]{Example}

\newtheorem{rem}[thm]{Remark}

\newtheorem*{ackno}{Acknowledgments} 

\theoremstyle{remark}

\newenvironment{say}{}{} 

\newenvironment{numerate}{%
  \begin{enumerate}}{%
    \end{enumerate}}

\author{Iku Nakamura} 
\address{Department of Mathematics
  \newline\indent
Hokkaido University
  \newline\indent
Sapporo 060-0810   \newline\indent
Japan}
\email{nakamura@math.hokudai.ac.jp}

\newcommand{\Prim}{\operatorname{Prim}}

\newcommand{\Rp}{_1R}

\newcommand{\wcL}{{\widetilde {\mathcal L}}}

\newcommand{\wtheta}{{\tilde\theta}}

\newcommand{\wTheta}{{\widetilde\Theta}}

\newcommand{\XR}{X_{\mathbb R}}

\newcommand{\Vor}{\operatorname{Vor}}
\newcommand{\Del}{\operatorname{Del}}

\newcommand{\St}{\operatorname{Star}}

\newcommand{\val}{\operatorname{val}}

\newcommand{\DDample}{\operatorname{DD_{ample}}}
\newcommand{\DDpol}{\operatorname{DD_{pol}}}
\newcommand{\DEGample}{\operatorname{DEG_{ample}}}
\newcommand{\DEGpol}{\operatorname{DEG_{pol}}}

\catcode`\@=11\relax
\newwrite\@unused
\def\typeout#1{{\let\protect\string\immediate\write\@unused{#1}}}
\def\figurepath{./}

\def\@nnil{\@nil}
\def\@empty{}
\def\@psdonoop#1\@@#2#3{}
\def\@psdo#1:=#2\do#3{\edef\@psdotmp{#2}\ifx\@psdotmp\@empty \else
    \expandafter\@psdoloop#2,\@nil,\@nil\@@#1{#3}\fi}
\def\@psdoloop#1,#2,#3\@@#4#5{\def#4{#1}\ifx #4\@nnil \else
       #5\def#4{#2}\ifx #4\@nnil \else#5\@ipsdoloop #3\@@#4{#5}\fi\fi}
\def\@ipsdoloop#1,#2\@@#3#4{\def#3{#1}\ifx #3\@nnil 
       \let\@nextwhile=\@psdonoop \else
      #4\relax\let\@nextwhile=\@ipsdoloop\fi\@nextwhile#2\@@#3{#4}}
\def\@tpsdo#1:=#2\do#3{\xdef\@psdotmp{#2}\ifx\@psdotmp\@empty \else
    \@tpsdoloop#2\@nil\@nil\@@#1{#3}\fi}
\def\@tpsdoloop#1#2\@@#3#4{\def#3{#1}\ifx #3\@nnil 
       \let\@nextwhile=\@psdonoop \else
      #4\relax\let\@nextwhile=\@tpsdoloop\fi\@nextwhile#2\@@#3{#4}}
\def\psdraft{
	\def\@psdraft{0}
}
\def\psfull{
	\def\@psdraft{100}
}
\psfull
\newif\if@prologfile
\newif\if@postlogfile
\newif\if@noisy
\def\pssilent{
	\@noisyfalse
}
\def\psnoisy{
	\@noisytrue
}
\psnoisy
\newif\if@bbllx
\newif\if@bblly
\newif\if@bburx
\newif\if@bbury
\newif\if@height
\newif\if@width
\newif\if@rheight
\newif\if@rwidth
\newif\if@clip
\newif\if@verbose
\def\@p@@sclip#1{\@cliptrue}

\def\@p@@sfile#1{\def\@p@sfile{null}%
	        \openin1=#1
		\ifeof1\closein1%
		       \openin1=\figurepath#1
			\ifeof1\typeout{Error, File #1 not found}
			\else\closein1
			    \edef\@p@sfile{\figurepath#1}%
                        \fi%
		 \else\closein1%
		       \def\@p@sfile{#1}%
		 \fi}
\def\@p@@sfigure#1{\def\@p@sfile{null}%
	        \openin1=#1
		\ifeof1\closein1%
		       \openin1=\figurepath#1
			\ifeof1\typeout{Error, File #1 not found}
			\else\closein1
			    \def\@p@sfile{\figurepath#1}%
                        \fi%
		 \else\closein1%
		       \def\@p@sfile{#1}%
		 \fi}

\def\@p@@sbbllx#1{
		\@bbllxtrue
		\dimen100=#1
		\edef\@p@sbbllx{\number\dimen100}
}
\def\@p@@sbblly#1{
		\@bbllytrue
		\dimen100=#1
		\edef\@p@sbblly{\number\dimen100}
}
\def\@p@@sbburx#1{
		\@bburxtrue
		\dimen100=#1
		\edef\@p@sbburx{\number\dimen100}
}
\def\@p@@sbbury#1{
		\@bburytrue
		\dimen100=#1
		\edef\@p@sbbury{\number\dimen100}
}
\def\@p@@sheight#1{
		\@heighttrue
		\dimen100=#1
   		\edef\@p@sheight{\number\dimen100}
}
\def\@p@@swidth#1{
		\@widthtrue
		\dimen100=#1
		\edef\@p@swidth{\number\dimen100}
}
\def\@p@@srheight#1{
		\@rheighttrue
		\dimen100=#1
		\edef\@p@srheight{\number\dimen100}
}
\def\@p@@srwidth#1{
		\@rwidthtrue
		\dimen100=#1
		\edef\@p@srwidth{\number\dimen100}
}
\def\@p@@ssilent#1{ 
		\@verbosefalse
}
\def\@p@@sprolog#1{\@prologfiletrue\def\@prologfileval{#1}}
\def\@p@@spostlog#1{\@postlogfiletrue\def\@postlogfileval{#1}}
\def\@cs@name#1{\csname #1\endcsname}
\def\@setparms#1=#2,{\@cs@name{@p@@s#1}{#2}}
\def\ps@init@parms{
		\@bbllxfalse \@bbllyfalse
		\@bburxfalse \@bburyfalse
		\@heightfalse \@widthfalse
		\@rheightfalse \@rwidthfalse
		\def\@p@sbbllx{}\def\@p@sbblly{}
		\def\@p@sbburx{}\def\@p@sbbury{}
		\def\@p@sheight{}\def\@p@swidth{}
		\def\@p@srheight{}\def\@p@srwidth{}
		\def\@p@sfile{}
		\def\@p@scost{10}
		\def\@sc{}
		\@prologfilefalse
		\@postlogfilefalse
		\@clipfalse
		\if@noisy
			\@verbosetrue
		\else
			\@verbosefalse
		\fi
}
\def\parse@ps@parms#1{
	 	\@psdo\@psfiga:=#1\do
		   {\expandafter\@setparms\@psfiga,}}
\newif\ifno@bb
\newif\ifnot@eof
\newread\ps@stream
\def\bb@missing{
	\if@verbose{
		\typeout{psfig: searching \@p@sfile \space  for bounding box}
	}\fi
	\openin\ps@stream=\@p@sfile
	\no@bbtrue
	\not@eoftrue
	\catcode`\%=12
	\loop
		\read\ps@stream to \line@in
		\global\toks200=\expandafter{\line@in}
		\ifeof\ps@stream \not@eoffalse \fi
		\@bbtest{\toks200}
		\if@bbmatch\not@eoffalse\expandafter\bb@cull\the\toks200\fi
	\ifnot@eof \repeat
	\catcode`\%=14
}	
\catcode`\%=12
\newif\if@bbmatch
\def\@bbtest#1{\expandafter\@a@\the#1
\long\def\@a@#1
\long\def\bb@cull#1 #2 #3 #4 #5 {
	\dimen100=#2 bp\edef\@p@sbbllx{\number\dimen100}
	\dimen100=#3 bp\edef\@p@sbblly{\number\dimen100}
	\dimen100=#4 bp\edef\@p@sbburx{\number\dimen100}
	\dimen100=#5 bp\edef\@p@sbbury{\number\dimen100}
	\no@bbfalse
}
\catcode`\%=14
\def\compute@bb{
		\no@bbfalse
		\if@bbllx \else \no@bbtrue \fi
		\if@bblly \else \no@bbtrue \fi
		\if@bburx \else \no@bbtrue \fi
		\if@bbury \else \no@bbtrue \fi
		\ifno@bb \bb@missing \fi
		\ifno@bb \typeout{FATAL ERROR: no bb supplied or found}
			\no-bb-error
		\fi
		\count203=\@p@sbburx
		\count204=\@p@sbbury
		\advance\count203 by -\@p@sbbllx
		\advance\count204 by -\@p@sbblly
		\edef\@bbw{\number\count203}
		\edef\@bbh{\number\count204}
}
\def\in@hundreds#1#2#3{\count240=#2 \count241=#3
		     \count100=\count240	
		     \divide\count100 by \count241
		     \count101=\count100
		     \multiply\count101 by \count241
		     \advance\count240 by -\count101
		     \multiply\count240 by 10
		     \count101=\count240	
		     \divide\count101 by \count241
		     \count102=\count101
		     \multiply\count102 by \count241
		     \advance\count240 by -\count102
		     \multiply\count240 by 10
		     \count102=\count240	
		     \divide\count102 by \count241
		     \count200=#1\count205=0
		     \count201=\count200
			\multiply\count201 by \count100
		 	\advance\count205 by \count201
		     \count201=\count200
			\divide\count201 by 10
			\multiply\count201 by \count101
			\advance\count205 by \count201
		     \count201=\count200
			\divide\count201 by 100
			\multiply\count201 by \count102
			\advance\count205 by \count201
		     \edef\@result{\number\count205}
}
\def\compute@wfromh{
		\in@hundreds{\@p@sheight}{\@bbw}{\@bbh}
		\edef\@p@swidth{\@result}
}
\def\compute@hfromw{
		\in@hundreds{\@p@swidth}{\@bbh}{\@bbw}
		\edef\@p@sheight{\@result}
}
\def\compute@handw{
		\if@height 
			\if@width
			\else
				\compute@wfromh
			\fi
		\else 
			\if@width
				\compute@hfromw
			\else
				\edef\@p@sheight{\@bbh}
				\edef\@p@swidth{\@bbw}
			\fi
		\fi
}
\def\compute@resv{
		\if@rheight \else \edef\@p@srheight{\@p@sheight} \fi
		\if@rwidth \else \edef\@p@srwidth{\@p@swidth} \fi
}
\def\compute@sizes{
	\compute@bb
	\compute@handw
	\compute@resv
}
\def\psfig#1{\vbox {
	\ps@init@parms
	\parse@ps@parms{#1}
	\compute@sizes
	\ifnum\@p@scost<\@psdraft{
		\if@verbose{
			\typeout{psfig: including \@p@sfile \space }
		}\fi
		\special{ps::[begin] 	\@p@swidth \space \@p@sheight \space
				\@p@sbbllx \space \@p@sbblly \space
				\@p@sbburx \space \@p@sbbury \space
				startTexFig \space }
		\if@clip{
			\if@verbose{
				\typeout{(clip)}
			}\fi
			\special{ps:: doclip \space }
		}\fi
		\if@prologfile
		    \special{ps: plotfile \@prologfileval \space } \fi
		\special{ps: plotfile \@p@sfile \space }
		\if@postlogfile
		    \special{ps: plotfile \@postlogfileval \space } \fi
		\special{ps::[end] endTexFig \space }
		\vbox to \@p@srheight true sp{
			\hbox to \@p@srwidth true sp{
				\hss
			}
		\vss
		}
	}\else{
		\vbox to \@p@srheight true sp{
		\vss
			\hbox to \@p@srwidth true sp{
				\hss
				\if@verbose{
					\@p@sfile
				}\fi
				\hss
			}
		\vss
		}
	}\fi
}}
\def\psglobal{\typeout{psfig: PSGLOBAL is OBSOLETE; use psprint -m instead}}
\catcode`\@=12\relax

\begin{document}
\bibliographystyle{amsalpha+}
\title[Mumford's construction]%
{On Mumford's construction of degenerating abelian varieties}



\begin{abstract}
  For a one-dimensional family of abelian varieties equipped with
  principal theta divisors a canonical limit is constructed as a pair
  consisting of a reduced projective variety and a Cartier divisor on
  it. Properties of such pairs are established.
\end{abstract}

\maketitle

\begin{center}
  {\small\bf August 19, 1997}
\end{center}

\tableofcontents


\section*{Introduction}
\label{sec:Introduction}

\begin{say}
  Assume that we are given a 1-parameter family of principally
  polarized abelian varieties with theta divisors. By this we will
  mean that we are in one of the following situations:
  \begin{enumerate}
  \item $\cR$ is a complete discrete valuation ring (DVR, for short)
    with the fraction field $K$, $S=\Spec \cR$, $\eta=\Spec K$ is the
    generic point, and we have an abelian variety $G_{\eta}$ over $K$
    together with an effective ample divisor $\Theta_{\eta}$ defining
    a principal polarization; or
  \item we have a projective family $(G,\Theta)$ over a small
    punctured disk $D_{\varepsilon}^0$.
  \end{enumerate}
  In this paper we show that, possibly after a finite ramified base
  change, the family can be completed in a simple and absolutely
  canonical manner to a projective family $(P,\Theta)$ with a
  relatively ample Cartier divisor $\Theta$ over $S$, resp.
  $D_{\varepsilon}$.  Moreover, this construction is stable under
  further finite base changes.  We give a combinatorial description of
  this family and its central fiber $(P_0,\Theta_0)$ and study their
  basic properties. In particular, we prove that $P_0$ is reduced and
  Cohen-Macaulay and that $H^i(P_0,\cO(d\Theta_0))$, $d\ge0$ are the
  same as for an ordinary PPAV (principally polarized abelian
  variety).
\end{say}

\begin{say}
  Existence of such construction has profound consequences for the
  moduli theory.  Indeed, with it one must expect that there exists a
  canonical compactification $\overline{A}_g$ of the moduli space
  $A_g$ of PPAVs, similar to the Mumford-Deligne compactification of
  the moduli space of curves. Without it, one has to believe that
  there is no single ``best'' geometrically meaningful
  compactification of $A_g$ and work with the infinitely many toroidal
  compactifications instead.  The moduli implications of our
  construction are explored in \cite{Alexeev_CMAV}.
\end{say}

\begin{say}
  Degenerations of abelian varieties have been studied exhaustively
  which makes our result somewhat surprising. There is a very complete
  description of degenerations of polarized abelian varieties of
  arbitrary degree of polarization over a complete Noetherian domain
  of arbitrary dimension, with or without an ample line bundle.  This
  description is called Mumford's construction, it was first published
  in a beautiful short paper \cite{Mumford_AnalyticDegsAVs} and later
  substantially expanded and improved by Faltings and Chai in
  \cite{Faltings85,Chai85,FaltingsChai90} (we note a parallel
  construction of Raynaud which works in the context of rigid analytic
  geometry). Mumford's construction gives an equivalence of categories
  $\DEGpol$, resp. $\DEGample$ of degenerations of polarized abelian
  varieties, resp. with a line bundle, and the categories $\DDpol$,
  resp. $\DDample$ of the ``degeneration data''.
\end{say}

\begin{say}
  As an auxiliary tool, Mumford's construction uses {\em relatively
    complete models.\/} Mumford remarks that such a model ``is neither
  unique nor canonical'' and that ``in fact, the non-uniqueness of
  $\wP$ gives one freedom to seek for the most elegant solutions in
  any particular case''. What we show in this paper is that if one is
  willing to give up some of the properties of $\wP$ and concentrate
  on the others, then in fact there is a canonical choice! Here is
  what we do:
  \begin{enumerate}
  \item We only consider the case of a 1-dimensional base $S$. This
    certainly makes the problem easier but not significantly.  In view
    of the moduli theory one shouldn't expect that a
    higher-dimensional family can be canonically completed, unless one
    is in a very special situation, such as for a ``test family'' over
    a special toric scheme.
  \item We allow an additional finite ramified base change $S'\to S$,
    even after one already has the semiabelian reduction.
    This, again, is perfectly natural from the moduli point of view.
  \item Most importantly, we do not care where in the central fiber
    the limit of the zero section of $G_{\eta}$ ends up. Our
    relatively complete model contains a semiabelian group scheme in
    many different ways, but the closure of the zero section need not
    be be contained in any of them. Hence, $P_0$ is a limit of
    $G_{\eta}$ as an abelian torsor, not as an abelian variety.
  \item Instead of a section, we pay a very special attention to the
    limit of the theta divisor $\Theta_{\eta}$, something which was
    overlooked in the previous constructions.
  \end{enumerate}
  For the most part of the paper we work in the algebraic situation,
  over a complete DVR. The complex-analytic case is entirely
  analogous, and we explain the differences in Section
  \ref{sec:Complex-analytic case}.
\end{say}

\begin{say}
  Shortly after \cite{Mumford_AnalyticDegsAVs} appeared, a series of
  works of Namikawa and Nakamura \cite{Namikawa_NewCompBoth,
    Namikawa_ToroidalDegsAVs, Namikawa_ToroidalDegsAVs2,
    Namikawa_ToroidalCompSiegel, Nakamura_ModuliSQAVs,
    Nakamura_CompNeronModel} was published that dealt with the
  complex-analytic situation. They contain a toric construction for an
  extended 1-parameter family. This construction is very similar to
  Mumford's, and the main difference is a substitute for the
  relatively complete model.  One unpleasant property of that
  substitute is that in dimension $g\ge5$ the central fiber need not
  be reduced.
\end{say}

\begin{say}
  When restricted to the complex-analytic setting, our construction
  has a lot in common with the Namikawa-Nakamura construction as well.
  The main difference is again the fact that we use and pay special
  attention to the theta divisor.  Our solution to the problems
  arising in dimension $g\ge5$ is simple -- a base change.  To
  underline the degree of dependence on the previous work, we call our
  construction {\em simplified Mumford's construction\/} and we call
  the central fibers {\em stable quasiabelian varieties,\/} or SQAVs,
  following Namikawa. We call the pairs $(P_0,\Theta_0)$ {\em stable
    quasiabelian pairs,\/} or SQAP. 
\end{say}

\begin{say}
  We note that Namikawa had constructed families $\cX^{(2n)}_g$ over
  the Voronoi compactification
  $\overline{A}^{\operatorname{VOR}}_{g,1,2n}$ of the moduli space
  $A_{g,1,2n}$ of PPAVs with a principal level structure of level
  $2n$, $n\ge3$. The boundary fibers in these families are different
  for different $n$, and some of them are non-reduced when $g\ge5$.
\end{say}

\begin{ackno}
  This research started when the first author was visiting the
  University of Tokyo and he would like to express his gratitude for
  the hospitality. We also thank Professors F. Bogomolov, Y.Kawamata,
  S.Mori and T.Oda for very stimulating conversations.
  
  The work of the first author was partially supported by the NSF. The
  work of the second author was partially supported by the
  Grant-in-Aid (No. 06452001) for Scientific Research, the Ministry of
  Education, Science, Sports and Culture, Japan.
\end{ackno}


\section{Delaunay and Voronoi polyhedral decompositions}
\label{sec:Delaunay and Voronoi polyhedral decompositions}

\begin{say}
  The structure of the extended family will be
  described explicitly in terms of two polyhedral decompositions which
  we now introduce.
\end{say}

\begin{notation}
  $X\simeq\bZ^r$ will denote a lattice in a real vector space
  $\XR=X\otimes\bR$, for a fixed positive integer $r$.  $B:X\times
  X\to\bR$ will be a symmetric bilinear form assumed to be positive
  definite. We denote the norm $\sqrt{B(x,x)}$ by $\|x\|_B$ or simply
  by $\|x\|$.
\end{notation}

\begin{defn}
  \label{defn:Delaunay cell}
  For an arbitrary $\alpha\in\XR$ we say that a lattice element $x\in
  X$ is $\alpha$-nearest if
  \begin{displaymath}
    \|x-\alpha\|_B= \min\{ \|x'-\alpha\|_B \quad\big|\, x'\in X \}.
  \end{displaymath}
  
  We define a $B$-Delaunay cell $\sigma$ (or simply a {\em Delaunay
    cell\/} if $B$ is understood) to be the closed convex hull of all
  lattice elements which are $\alpha$-nearest for some fixed
  $\alpha\in\XR$.  Note that for a given Delaunay cell $\sigma$ the
  element $\alpha$ is uniquely defined only if $\sigma$ has the
  maximal possible dimension, equal to $r$. In this case $\alpha$ is
  called the {\em hole \/} of $\sigma$, cf. Section 2.1.2 of the
  ``encyclopedia of sphere packings and lattices''
  \cite{ConwaySloane93}.  One should imagine a sphere around the
  $\alpha$-closest lattice elements (which is known as ``the empty
  sphere'' because there are no other lattice elements in its
  interior) with $\alpha$ at the center.
  
  Together all the Delaunay cells constitute a locally finite
  decomposition of $\XR$ into infinitely many bounded convex polytopes
  which we call the {\em Delaunay cell decomposition \/} $\Del_B$.
\end{defn}

\begin{rem}
  It is clear from the definition that the Delaunay decomposition is
  invariant under translation by the lattice $X$ and that the
  0-dimensional cells are precisely the elements of $X$.
\end{rem}

\begin{defn}
  For a given $B$-Delaunay cell $\sigma$ consider all $\alpha\in\XR$
  that define $\sigma$. They themselves form a locally closed bounded
  convex polytope.  We denote the closure of this polytope by
  $\hat\sigma=V(\sigma)$ and call it the $B$-Voronoi cell or simply
  the {\em Voronoi cell\/}. The Voronoi cells make up the {\em Voronoi
    cell decomposition $\Vor_B$\/} of $\XR$.
\end{defn}

\begin{say}
  It is easy to see that the Delaunay and Voronoi cells are dual to
  each other in the following sense:
\end{say}

\begin{lem}
  \label{lem:Voronoi_Delaunay_duality}
  \begin{numerate}
  \item For a fixed form $B$ there is a 1-to-1 correspondence between
    Delaunay and Voronoi cells given by $\hat\sigma=V(\sigma)$,
    $\sigma=D(\hat\sigma)$.
  \item $\dim\sigma+\dim\hat\sigma=r$.
  \item $\sigma\subset\tau$ if and only if $\hat\tau\subset\hat\sigma$.
  \item For $\sigma=x\in X$ the corresponding Voronoi cell $V(x)$ has
    the maximal dimension. $V(x)$ is the set of points of $\XR$ that
    are at least as close to $x$ as to any other lattice element $x'$.
  \item For an arbitrary Delaunay cell $\sigma$ the dual Voronoi cell
    $\hat\sigma$ is the polytope with vertices at holes $\alpha(\tau)$
    where $\tau$ goes over all maximal-dimensional Delaunay cells
    containing $\sigma$.
  \end{numerate}
\end{lem}

\begin{rem}
  With the definition this natural, it is no wonder that Voronoi and
  Delaunay cells have a myriad of applications in physics, chemistry
  and even geography (\cite{OkabeBootsSugihara92}, many more
  references in \cite{ConwaySloane93}) and go by many different names.
  Some alternative names for Voronoi cells are: Voronoi polytopes,
  nearest neighbor regions, Dirichlet regions, Brillouin zones,
  Wigner-Seitz cells. Delaunay cells had been called by various
  authors Delone cells, Delony cells, $L$-polytopes. Evidently, even
  the spelling of the last name of Boris Nikolaevich Delone is not
  agreed upon.  A quick computer database search shows that the French
  variant ``Delaunay'' is prefered in 99\% of all papers, so this is
  our choice too.
\end{rem}

\begin{exmp}
  The Figures 1 and 2 give the only two, up to the action of
  $SL(2,\bZ)$, Delaunay decompositions of $\bZ^2$. The broken lines
  show the corresponding Voronoi decompositions. Note that, unlike
  Delaunay, the Voronoi decompositions may have some continuous
  moduli.
\end{exmp}

 
\begin{say}
  Here is another way to understand the Delaunay decompositions.
\end{say}

\begin{lem}
  \label{lem:paraboloid}
  Consider an $(r+1)$-dimensional real vector space $\XR\oplus\bR$
  with coordinates $(x,x_0)$ and a paraboloid in it defined by the
  equation
  \begin{displaymath}
    x_0=A(x)=B(x,x)/2+lx/2 
  \end{displaymath}
  for some $l\in\XR^*$. Consider the convex hull $Q$ of countably many
  points on this paraboloid with $x\in X$.
  
  This object has a multifaceted shape and projections of the facets
  onto $\XR$ are precisely the $B$-Delaunay cells. The equations that
  cut out the cone at the vertex $(c,A(c))$ are
  \begin{displaymath}
    x_0-A(c)\ge B\big( \alpha(\sigma),x\big) +lx/2=
    dA(\alpha(\sigma))(x) , 
  \end{displaymath}
  where $\sigma$ goes over all the maximal-dimensional Delaunay cells
  containing $c$. If $\sigma$ is a maximal-dimensional Delaunay cell,
  then the interior normal in the dual space $\XR^*\oplus\bR$ to the
  corresponding facet is $\big(1,-dA(\alpha(\sigma))\big)$. The normal
  fan to the paraboloid consists of $\{0\}$ and the cones over the
  shifted Voronoi decomposition $\big(1,-dA(\Vor_B)\big)$.
\end{lem}

 
\begin{say}
  The inequality for $x_0$ does not depend on the cell $\sigma\ni x$
  chosen because for any Delaunay vector $v\in\sigma_1\cap\sigma_2$
  \begin{displaymath}
    B(\alpha(\sigma_1),v)=B(v,v)/2=B(\alpha(\sigma_2),v)
  \end{displaymath}
  (see \ref{rem:equations_for_hole}). We look at $dA$ here as being a
  map from $\XR$ to $\XR^*$.
\end{say}

\begin{proof}
  We claim that the hyperplanes that cut out the facets at the origin
  are in 1-to-1 correspondence with the maximal-dimensional Delaunay
  cells containing 0. Let $\sigma\ni0$ be one of such cells with a
  hole $\alpha$. Then the points $(x,A(x))$ with $x\in\sigma$ lie on
  the hyperplane $x_0=B(\alpha,x)+lx/2$ and those with $x\notin\sigma$
  lie above it. Indeed,
  \begin{displaymath}
    A(x)-B(\alpha,x)-lx/2=
    B(x,x)/2-B(\alpha,x)=(\|\alpha-x\|^2-\|\alpha-0\|^2)/2 \ge0
  \end{displaymath}
  and the equality holds if and only if $x\in\sigma$. The other
  vertices are checked similarly.
\end{proof}

\begin{rem}\label{rem:equations_for_hole}
  As the proof shows, for a maximal-dimensional Delaunay cell the hole
  $\alpha$ is the unique solution of the system of linear equations
  \begin{displaymath}
    2B(\alpha,x)=B(x,x), \quad x\in\sigma\cap X.
  \end{displaymath}
\end{rem}

\begin{defn}\label{defn:eta_xc}
  We will denote the minimum of the functions $dA(\alpha(\sigma))(x)$
  in the above lemma by $\eta(x,c)$. 
\end{defn}

\begin{rem}
  Note that all $\eta(x,c)$ are integral-valued on $X$ if and only if
  all $\alpha(\sigma)$ are integral, i.e., belong to $X^*$.
\end{rem}

\begin{defn}\label{defn:star_Delaunay_primitive}
  The union of the Delaunay cells containing a lattice vertex $c$ is
  called the {\em star of Delaunay cells\/} and denoted by $\St(c)$.
  For a cell $\sigma\subset\St(0)$ the nonzero lattice elements $v_1\dots
  v_m\in\sigma$ are called the Delaunay vectors. In addition, we
  introduce
  \begin{displaymath}
    \Cone(0,\sigma)=\bR_+v_1+\dots+\bR_+v_m,
  \end{displaymath}
  a convex cone with the vertex at the origin. By translation we
  obtain cones $\Cone(c,\sigma)$ with vertices at other lattice elements.
  The function $\eta(x,c)$ is linear on each $\Cone(c,\sigma)$,
  $\sigma\subset\St(c)$.

  A lattice vector in a cone is called {\em primitive\/} if it cannot
  be written as a sum of two nonzero lattice vectors in this cone. We
  will denote by $\Prim$ the union of primitive vectors in Delaunay
  cells $\sigma$ for all $\sigma\subset\St(0)$.

  We say that some elements $x_1\dots x_m$ are {\em cellmates\/} if
  there exists a (maximal-dimensional) Delaunay cell $\tau\ni0$ such
  that $x_1\dots x_m\in \Cone(0,\tau)$.
\end{defn}
  
\begin{say}
  We thank S. Zucker for suggesting the term ``cellmates''.
\end{say}

\begin{defn}\label{defn:generating_cell}
  We call the cell $\sigma$ {\em generating\/} if its Delaunay
  vectors contain a basis of $X$ and {\em totally generating\/} if,
  moreover, the integral combinations $\sum n_iv_i$ of Delaunay
  vectors with $n_i\ge0$ give all lattice elements of the cone
  $\Cone(0,\sigma)$.
  
\end{defn}

\begin{saynum}
  In dimensions $g\le4$ all Delaunay cells are totally generating as
  was proved by Voronoi in a classical series of papers
  \cite{Voronoi08all}. It was only recently shown that there exist
  non-generating cells in dimension 5 (see \cite{ErdahlRyshkov87},
  p.796). But the following example from \cite{Erdahl92} of a
  non-generating cell is by far the easiest.
\end{saynum}

\begin{exmp}\label{exmp:E_8_not_generating}
  Take $B=E_8$, an even unimodular positive definite quadratic form
  given by a familiar $8\times8$ matrix with +2 on the main diagonal.
  By \ref{rem:equations_for_hole} we have
  \begin{displaymath}
    E_8(\alpha,v_i)=E_8(v_i,v_i)/2\in\bZ, \quad v_i\in\sigma.
  \end{displaymath}
  If the Delaunay vectors $v_i$ contained a basis then we would have
  $\alpha\in X$, since $E_8$ is unimodular. But the holes cannot
  belong to the lattice by their very definition.
  
  In this example each maximal-dimensional Delaunay cell generates a
  sublattice of index 2 or 3.
\end{exmp}

\begin{defn}\label{defn:nilpotency}
  The {\em nilpotency\/} of a Delaunay cell $\sigma$ is the minimal
  positive integer $n$ such that the lattice generated by the vectors
  $v_1/n,\dots, v_m/n$ contains $X$. The nilpotency of the Delaunay
  decomposition is the least common multiple of nilpotencies of its
  cells.
\end{defn}

\begin{rem}
  Existence of non-generating cells is responsible for many
  combinatorial complications that occur in dimension $g\ge5$.
\end{rem}



\section{Degeneration data}
\label{sec:Degeneration data}

\begin{say}
  The purpose of this section is to recall the Mumford-Faltings-Chai
  description \cite[III]{FaltingsChai90} of the degenerations of
  abelian varieties. We will use this description in the next section.
  We refer the reader to \cite{FaltingsChai90} for basic definitions
  and facts about polarizations, semiabelian group schemes etc.  For
  easier reference we use the same notation as in
  \cite{FaltingsChai90} where this is convenient.
\end{say}

\begin{notation}
  $\cR$ is a Noetherian normal integral domain complete with respect
  to an ideal $I=\sqrt{I}$, $K$ is the fraction field and $k=\cR/I$ is
  the residue ring, $S=\Spec \cR$, $S_0=\Spec k$, and $\eta=\Spec K$
  is the generic point.  $\cR$ is assumed to be regular or to be the
  completion of a normal excellent domain.
  
  In our application $\cR$ will be a DVR, and $k$ will be the residue
  field. 
\end{notation}

\begin{defn}
  The objects of the category $\DEGample$ are pairs $(G,\cL)$, where
  $G$ is a semiabelian group scheme over $S$ with abelian $G_{\eta}$
  and $\cL$ is an invertible sheaf on $G$ with ample $\cL_{\eta}$. The
  morphisms are group homomorphisms respecting $\cL$'s.
  
  The objects of the category $\DEGpol$ are pairs $(G,\lambda_{\eta})$
  with $G$ as above and with the polarization
  $\lambda_{\eta}:G_{\eta}\to G^t_{\eta}$, where $G^t_{\eta}$ is the
  dual abelian variety. The morphisms are group homomorphisms
  respecting $\lambda$'s.
\end{defn}

\begin{say}
  Next, we recall the definition of the degeneration data. We will
  only need the split case.
\end{say}

\begin{defn}
  The degeneration data in the split case consists of the following:
  \begin{enumerate}
  \item An abelian scheme $A/S$ of relative dimension $a$, a split
    torus $T/S$ of relative dimension $r$, $g=a+r$, and a semiabelian
    group scheme $\wG/S$,
    \begin{displaymath}
      1\to T\to \wG \overset{\pi}{\to} A \to 0.
    \end{displaymath}
    This extension is equivalent via negative of pushout to a
    homomorphism $c:\uX\to A^t$.  Here $X$ is a rank $r$ free
    commutative group, and $\uX=X_S$ is a constant group scheme, the
    group of characters of $T$. $A^t/S$ is the dual abelian scheme
    of~$A/S$.
  \item A rank $r$ free commutative group $Y$ and the constant group
    scheme $\uY=Y_S$.
  \item A homomorphism $c^t:\uY\to A$. This is equivalent to giving an
    extension
    \begin{displaymath}
      1\to T^t\to \wG^t \to A^t \to 0,
    \end{displaymath}
    where $T^t$ is a torus with the group of characters $\uY$.
  \item An injective homomorphism $\phi:Y\to X$ with finite cokernel.
  \item A homomorphism $\iota: Y_{\eta}\to\wG_{\eta}$ lying over
    $c^t_{\eta}$. This is equivalent to giving a bilinear section
    $\tau$ of $(c^t\times c)^* \cP^{-1}_{A,\eta}$ on $Y\times X$, in
    other words a trivialization of the biextension $\tau:1_{Y\times
      X}\to (c^t\times c)^*\cP^{-1}_{A,\eta}$. Here, $\cP_{\eta}$ is
    the Poincare sheaf on $A_{\eta}\times A^t_{\eta}$, which comes
    with a canonical biextension structure.
  \item An ample sheaf $\cM$ on $A$ inducing a polarization
    $\lambda_A:A\to A^t$ of $A/S$ such that $\lambda_Ac^t=c\phi$. This
    is equivalent to giving a $T$-linearized sheaf $\wcL=\pi^*\cM$ on
    $\wG$. 
  \item An action of $Y$ on $\wcL_{\eta}$ compatible with $\phi$. This
    is equivalent to a cubical section $\psi$ of
    $(c^t)^*\cM_{\eta}^{-1}$ on $Y$, in other words to a cubical
    trivialization $\psi:1_{Y}\to(c^t)^*\cM_{\eta}^{-1}$, which is
    compatible with $\tau\circ(\id_Y\times\phi)$. $\psi$ is defined up
    to a shift by $Y$.
  \end{enumerate} 
  The trivialization $\tau$ is required to satisfy the following {\em
    positivity condition}: $\tau(y,\phi y)$ for all $y$
  extends to a section of $\cP^{-1}$ on $A\underset{S}{\times} A^t$,
  and it is 0 modulo $I$ if $y\ne0$.
  
  The objects of the category $\DDample$ are the data above, and the
  morphisms are the homomorphisms of $\wG$'s respecting this data.
\end{defn}

\begin{defn}
  Similarly, the objects of the category $\DDpol$ consist of the data
  as above minus the sheaves $\cM$, $\wcL$ and the section $\psi$,
  with the positivity condition again. In addition, one requires the
  trivialization $\tau$ to be symmetric (in the previous case this was
  automatic). The morphisms are homomorphisms of $\wG$'s respecting
  this data.
\end{defn}

\begin{thm}[Faltings-Chai]\label{thm:Faltings_Chai_equivalence}
  The categories $\DEGample$ and $\DDample$, resp. $\DEGpol$ and
  $\DDpol$ are equivalent.
\end{thm}

\begin{say}
  In the case where $A=0$ and $\wG=T$ is a torus (or, more generally,
  when $c=c^t=0$) the section $\tau$ is simply a bilinear function
  $b:Y\times X\to K^*$, and $\psi$ is a function $a:Y\to K^*$
  satisfying
  \begin{enumerate}
  \item $b(y_1,\phi y_2)=a(y_1+y_2)a(y_1)^{-1}a(y_2)^{-1}$,
  \item $b(y,\phi y)\in I$ for $y\ne0$.
  \end{enumerate}
  This also implies that
  \begin{enumerate}\setcounter{enumi}{2}
  \item $b(y_1,\phi y_2) = b(y_2, \phi y_1)$,
  \item $a(0)=1$,
  \item $a(y_1+y_2+y_3)
    a(y_2+y_3)^{-1} a(y_3+y_1)^{-1} a(y_1+y_2)^{-1}
    a(y_1) a(y_2) a(y_3) =1$.
  \end{enumerate}
\end{say}

\begin{say}
  We will use the following notation.  $\wG$ is affine over $A$ and
  one has $\wG = \Spec_A (\oplus_{x\in X} \cO_x)$. Each $\cO_x$ is an
  invertible sheaf on $A$, canonically rigidified along the zero
  section, and one has $\cO_x\simeq c(x)$.  The pushout of the
  $T$-torsor $\wG$ over $A$ by $x\in X$ is $\cO_{-x}$.  The sheaf
  $\cM\otimes\cO_x$, rigidified along the zero section, is denoted by
  $\cM_x$.  The $Y$-action on $(\oplus_{x\in
    X}\cM_x)\underset{\cR}{\otimes}K$ defined by $\psi$ is denoted by
  $S_y:T^*_{c^t(y)}\cM_x \to \cM_{x+\phi y,\eta}$.
\end{say}


\section{Simplified Mumford's construction}
\label{sec:Simplified Mumford's construction}

\begin{setup}
  In this section, $\cR$ is a DVR, $I=(s)$, and $k$ is the residue
  field. We will denote the point $S_0$ simply by $0$. We start with
  an abelian variety $A_{\eta}$ with an effective ample Cartier
  divisor $\Theta_{\eta}$ defining principal polarization, and
  $\cL_{\eta}=\cO(\Theta_{\eta})$. Applying the stable reduction
  theorem (\cite{SGA71,ArtinWinters_StableReduction}) after a finite
  base change $S'\to S$ we have a semiabelian group scheme $G'/S'$ and
  an invertible sheaf $\cL$ extending $(A_{\eta}',\cL_{\eta}')$ such
  that the toric part $T'_0$ of the central fiber $G'_0$ is split. In
  order not to crowd notation, we will continue to denote the objects
  by $S$, $G$, $\cL$ etc.
  
  We have an object of $\DEGample$, and, by the previous section, an
  object of $\DDample$, i.e., the degeneration data. Since the
  polarization is principal, $\phi:Y\to X$ is an isomorphism and we
  can identify $Y$ with $X$. Further, the sheaf $\cM$ on $A$ defines a
  principal polarization. We denote by $\theta_A$ a generator of
  $H^0(A,\cM)$.
\end{setup}

\begin{say}
  Here is the main object of our study:
\end{say}

\begin{defn}
  Consider the graded algebra 
  \begin{displaymath}
    \cS_2=\big(\sum_{d\ge0}(\oplus_{x\in
      X}\cO_x)\otimes \cM^d\vartheta^d\big)\underset{R}{\otimes}K,    
  \end{displaymath}
  where $\vartheta$ is an indeterminate defining the grading. In this
  algebra consider the subalgebra $\Rp$ generated in degree by the
  $\cM_0=\cM$ and all its $Y$-translates, $S^*_y(\cM_0)$. This is a
  locally free graded $\cO_A$-algebra.  Finally, the algebra $R$ is
  the saturation of $\Rp$ in an obvious sense which will be further
  explained below. We define the scheme $\wP=\Proj_A R$ and the sheaf
  $\wcL$ on it as $\cO(1)$.
  
  For each $x\in X=Y$ we have an element $S_x^*(\theta_A)\in
  H^0(A,\cM_x)$ that will be denoted by $\xi_x$. We have a formal
  power series
  \begin{displaymath}
    \wtheta= \sum_{x\in X} \xi_x.
  \end{displaymath}
\end{defn}

\begin{say}
  We will see that, possibly after another finite base change $S'\to
  S$, the scheme $\wP$ is a {\em relatively complete model\/} as
  defined in \cite[III.3]{FaltingsChai90}. Via Mumford's construction,
  this gives a projective scheme $P/S$ extending $A_{\eta}$. We will
  see that it naturally comes with a relative Cartier divisor
  $\Theta$.
\end{say}

\begin{say}
  The subalgebra $R$ defines a subalgebra $R'$ in 
  \begin{displaymath}
    \cS_1=\big(\sum_{d\ge0}(\oplus_{x\in
      X}\cO_x)\vartheta^d\big)\underset{R}{\otimes}K.    
  \end{displaymath}
  One has $\Proj R=\Proj R'$ and $\wcL\simeq\wcL'\otimes\pi^*\cM$.
\end{say}

\subsection{Case of maximal degeneration}
\label{subsec:simplified_max_deg}
\skipline

\begin{say}
  In this case, $A=0$ and $\wG=T=\Spec R[w^x \,;\, x\in X]$. $\Rp$ is
  the $R$-subalgebra of $K[\vartheta,w^x \,;\, x\in X]$ generated by
  $a(x)\vartheta$. We have 
  \begin{displaymath}
    \wtheta=\sum_{x\in X}\xi_x
    = \sum_{x\in X} a(x)w^x \vartheta.
  \end{displaymath}
\end{say}

\begin{defn}
  We define the functions $a_0:Y\to k^*$, $A:Y\to\bZ$, $b_0:Y\times
  X\to k^*$, $B:Y\times X\to\bZ$ by setting
  \begin{displaymath}
    a(y)=a'(y)s^{A(y)}, \qquad b(y,x)=b'(y,x)s^{B(y,x)}
  \end{displaymath}
  with $a'(y),b'(y) \in R\setminus I$ and taking $a_0,b_0$ to be
  $a',b'$ modulo $I$. 
\end{defn}

\begin{rem}
  We are using the letter $A$ for two purposes now: to denote an
  abelian variety, and to denote the integral-valued function
  above. This should not lead to any confusion since their meanings
  are very different.
\end{rem}

\begin{say}
  Through our identification $\phi:Y\isoto X$ the functions $a,A,b$
  and $B$ become functions on $X$ and $X\times X$.  The functions $a$
  and $A$ are quadratic non-homogeneous, the functions $b,B$ are
  symmetric and they are the homogeneous parts of $a^2,2A$,
  respectively. We have
  \begin{displaymath}
    A(x)=B(x,x)/2+lx/2
  \end{displaymath}
  for some $l\in X^*$. The positivity condition implies that $B$ is
  positive definite.
\end{say}

\begin{rem}
  The function $B:Y\times X\to\bZ$ describes the monodromy of the
  family $G$ and is called the monodromy pairing, cf.
  \cite[IX.10.4]{SGA71}.
\end{rem}

\begin{say}
  Since all $a'(y)$ are invertible in $\cR$, the algebra $\Rp$ is
  generated by monomials $\zeta_x=s^{A(x)}w^x\vartheta$, so it is a
  semigroup algebra.
\end{say}

\begin{defn}
  We introduce two lattices $M=X\oplus\bZ e_0\simeq \bZ^{r+1}$ and its
  dual $N=X^*\oplus\bZ f_0$.
\end{defn}

\begin{saynum}
  Each $\zeta_x$ corresponds to a lattice element $(x,A(x))\in M$.
  These are exactly the vertices of the multifaceted paraboloid $Q$ in
  Figure~3 
  which we imagine lying in the hyperplane $(1,M)$ inside $\bZ\oplus
  M$. The extra $\bZ$ corresponds to the grading by $\vartheta$. The
  saturation $R$ of $\Rp$ is generated by monomials corresponding to
  all lattice vectors lying inside $\Cone(Q)$.
\end{saynum}

\begin{thm}\label{thm:first_structure_thm_family}
  \begin{numerate}
  \item $\wP$ is covered by the affine toric schemes $U(c)=\Spec
    R(c)$, $c\in X$, where $R(c)$ is the semigroup algebra
    corresponding to the cone at the vertex $c\in Q$ of lattice
    elements
    \begin{displaymath}
      \{(x,x_0) \,|\, x_0 \ge \eta(x,c) \}
    \end{displaymath}
    ($\eta(x,c)$ is defined in \ref{defn:eta_xc}).
  \item $R(c)$ is a free $\cR$-module with the basis
    $\zeta_{x,c}=s^{\ulcorner\eta(x,c)\urcorner}w^x$ (here $\ulcorner
    z\urcorner$ denotes the least integer $\ge z$).
    \label{numi_basis_Rc}
  \item All the rings $R(c)$ are isomorphic to each other, and each is
    finitely generated over $\cR$. The scheme $\wP$ is locally of
    finite type over $\cR$.
  \item $\Spec R(c)$ is the affine torus embedding over $S=\Spec \cR$
    corresponding to the cone $\Delta(c)$ over
    \begin{displaymath}
      \big(1,-dA(\hc) \big) \subset(1,N_0) \subset N, 
    \end{displaymath}
    where $\hc$ is the Voronoi cell dual to $c$.
  \item $\wP$ is the torus embedding $T_N\emb\Delta$ where $\Delta$ is
    the fan in $N_{\bR}$ consisting of $\{0\}$ and the cones over the
    shifted Voronoi decomposition $(1,-dA(\Vor_B))$. The morphism
    $\wP\to S$ is described by the map of fans from $\Delta$ to the
    half line $\bR_{\ge0}f_0$.
  \item $\wcL$ is invertible and ample.
  \item One has natural compatible actions of $T$ on $\wP$ and of
    $T\times\bG_m$ on $\wcL$.
  \end{numerate}
\end{thm}

\begin{rem}
  The reference for torus embeddings over a DVR is \cite[IV,\S
  3]{ToroidalEmbeddingsI}. Formally, all computations work the same
  way as for a toric variety with torus action of $k[s,1/s,w^x \,;\,
  x\in X]$. 
\end{rem}

\begin{proof}
  The first part of (i) is simply the description of the standard
  cover of $\Proj$ by $\Spec$'s. The second part, as well as (iv) and
  (v) follow immediately from lemma \ref{lem:paraboloid}. (ii) and
  (iii) follow at once from (i).
  
  The ring extension $\Rp\subset R$ is integral, hence $\Proj R\to
  \Proj\, \Rp$ is well-defined and is finite. The sheaf $\cO(1)$ on
  $\Proj\, \Rp$ is invertible and ample, and $\wcL$ is its pullback.
  This gives (vi).
  
  The actions in (vii) are defined by the $X$-, resp.
  $(X\oplus\bZ)$-gradings.
\end{proof}

\begin{say}
  As a consequence, we can apply the standard description in the
  theory of torus embeddings of the open cover, torus orbits and their
  closures:
\end{say}

\begin{thm}\label{thm:toric_descr_family}
  \begin{numerate}
  \item For each Delaunay cell $\sigma\in\Del_B$ one has a ring
    $R(\sigma)$ corresponding to the cone over the dual Voronoi cell
    $\hat\sigma$.  $U(\sigma)$ is open in $\wP$ and $U(\sigma_1)\cap
    U(\sigma_2)= U(D(\hat\sigma_1\cap \hat\sigma_2))$.
  \item $R(\sigma)$ is the localization of $R(c)$, $c\in\sigma$, at
    $\zeta_{d,c}$, $d\in X\cap \bR(\sigma-c)$.
  \item In the central fiber $\wP_0$, the $T_0$-orbits are in 1-to-1
    dimension-preserving correspondence with the Delaunay cells
    $\sigma$. In particular, the irreducible components of $\wP$
    correspond to the maximal-dimensional cells.
  \item The closure $\wV(\sigma)$ of $\orb(\sigma)$ together with the
    restriction of the line bundle $\wcL$ is a projective toric variety
    over $k$ with a $T_0$-linearized ample line bundle corresponding
    to the lattice polytope $\sigma$.
  \item $\wV(\sigma_1)\cap \wV(\sigma_2)= \wV(\sigma_1\cap \sigma_2)$.
  \item For a maximal-dimensional cell $\sigma$, the multiplicity of
    $\wV(\sigma)$ in $\wP_0$ is the denominator of
    $dA(\alpha(\sigma))\in X^*_{\bQ}$.
  \end{numerate}
\end{thm}

\begin{question}
  For a maximal-dimensional cell $\sigma$, when is $\wP$ generically
  reduced at $\wV(\sigma)$? In other words, when is $dA(\alpha(\sigma))$
  integral?
\end{question}

\begin{lem}
  $dA(\alpha(\sigma))\in X^*$ in any of the following cases:
  \begin{numerate}
  \item $\sigma$ is generating.
  \item $A(x)/n\in\bZ$ for every $x\in X$, where $n$ is the
    nilpotency of $\sigma$.
  \end{numerate}
\end{lem}
\begin{proof}
  (i) is a particular case of (ii), so let us prove the second part.
  
  $dA(\alpha)\in X^*$ if and only if $dA(\alpha)(x)\in\bZ$ for every
  $x\in X$.  Now let $v_1\dots v_m$ be the Delaunay vectors of
  $\sigma$. By the definition of the nilpotency in
  Definition~\ref{defn:nilpotency} we have $x=(1/n) \sum n_iv_i$ for
  some $n_i\in\bZ$. Then
  \begin{align*}
    dA(\alpha)(x) &= B(\alpha,x) + lx/2 = 
    \frac{1}{n}\sum n_i\big( B(\alpha,v_i) +lv_i/2 \big)\\
    &= \frac{1}{n}\sum n_i \big( B(v_i,v_i)/2 +lv_i/2 \big) 
    = \sum n_i A(v_i)/n \in\bZ  
  \end{align*}
\end{proof}

\begin{saynum}
  For the central fiber $\wP_0$ to be generically reduced, we need
  $A(x)$ to be divisible by the nilpotency of the lattice. This
  certainly holds after a totally ramified base change.  Consider the
  polynomial $z^n-s\in K[z]$. It is irreducible by the Eisenstein
  criterion. The field extension $K\subset K'=K[z]/(z^n-s)$ has degree
  $n$ and is totally ramified. The integral closure of $R$ in $K'$ is
  again a DVR, complete with respect to the maximal ideal $I'=(s')$
  (see e.g.  \cite[II.3]{Serre79}) and $\val_{s'}(s)=n$, so that
  $A'(x)=nA(x)$.

  The following example, very similar to Example
  \ref{exmp:E_8_not_generating}, shows that this base change is indeed
  sometimes necessary.
\end{saynum}

\begin{exmp}
  Consider the degeneration data $A(x)=E_8(x)/2$, $B(y,x)=E_8(y,x)$.
  Then $dA=E_8$ and for every hole $\alpha(\sigma)$ we have
  $dA(\alpha)\notin X^*$. Indeed, otherwise we would have $\alpha\in
  X$ since $E_8$ is unimodular, and this is impossible by the
  definition of a hole.
  
  In this example every irreducible component of the central fiber
  $\wP_0$ has multiplicity $2$ or $3$.
\end{exmp}

\begin{assume}\label{assume:base_change_made}
  From now on, we assume that the necessary base change has been done,
  so $dA(\alpha)$ is integral for each hole $\alpha$. This implies
  that all $\eta(x,c)$ are integral-valued on $X$.
\end{assume}

\begin{saynum}\label{saynum:Yaction_wP}
  {\bf $Y$-action on $\wP$.} We are given a canonical $Y$-action on
  $K[\vartheta,w^x \,;\,x\in X]$ by construction. It is constant on
  $K$ and sends each generator $\xi_x=a(x)w^x\vartheta$ to another
  generator $\xi_{x+y}=a(x+y)w^{x+y}\vartheta$. Precisely because
  $a(x)$ is quadratic, this action extends uniquely to the whole
  $K[\vartheta,w^x]$.  Clearly, the subrings $\Rp$ and $R$ are
  $Y$-invariant, so we have the $Y$-action on $R$ which will be
  denoted by $S^*_y$. We easily compute:
\end{saynum}

\begin{eqnarray*}
  && S^*_y( \frac{a'(x+c)}{a'(x)}\zeta_{x,c} ) = 
  \frac{a'(x+c+y)}{ a'(x+y)} \zeta_{x,c+y},
  \\
  && S^*_y(\zeta_{x,c} ) = b'(x,y) \zeta_{x+y,c}.
\end{eqnarray*}

\begin{say}
  This describes the action $S_y^*:R(c)\to R(c+y)$ and $S_y:\Spec
  R(c+y) \to \Spec R(c)$.
\end{say}

\begin{thm}\label{thm:structure_of_central_fiber}
  $\wP_0$ is a scheme locally of finite type over $k$. It is covered
  by the affine schemes $\Spec R_0(\sigma)$ of finite type over $k$
  for Delaunay cells $\sigma\in\Del_B$, where
  $R_0(\sigma)=R(\sigma)\underset{\cR}{\otimes}k$. The following
  holds:
  \begin{numerate}
  \item $R_0(c)$ is a $k$-vector space with basis
    $\{\bar\zeta_{x,c}\,;\,x\in X\}$, the multiplication being defined
    by
    \begin{displaymath}
      \bar\zeta_{x_1,c}\dots\bar\zeta_{x_m,c}=\bar\zeta_{x_1+\dots+x_m,c}      
    \end{displaymath}
    if $x_1\dots x_m$ are cellmates with respect to the $B$-Delaunay
    decomposition, and $0$ otherwise.
    \label{enumi:basis_R0}
  \item For $\sigma\ni c$ the ring $R_0(\sigma)$ is the localization
    of $R_0(c)$ at $\{\bar\zeta_{d,c}\,;\, d\in X\cap\bR(\sigma-c)\}$.
  \item $U_0(\sigma_1)\cap U_0(\sigma_2)=
    U_0(D(\hat\sigma_1\cap\hat\sigma_2))$.
  \item The group $Y$ of periods acts on $\wP_0$ by sending $\Spec
    R_0(c+\phi(y))$ to $\Spec R_0(c)$ in the following way:
    \begin{displaymath}
      S^*_y(\bar\zeta_{x,c})= 
      b_0(y,x)\, \bar\zeta_{x,c+\phi(y)}.
    \end{displaymath}
  \end{numerate}
\end{thm}
\begin{proof}
  This follows from \ref{thm:first_structure_thm_family},
  \ref{thm:toric_descr_family} and \ref{saynum:Yaction_wP}.
  
  Here is a way to see (i) geometrically: each $\zeta_{x_i,c}$
  corresponds to a point on a face of the cone of $Q$ at $c$. The sum
  of several such points lie on a face if and only if they belong to a
  common face, i.e., if and only if $x_1\dots x_m$ are cellmates.
  Otherwise, the product corresponds to a point in the interior of the
  cone and equals $s^n\zeta_{x_1+\dots+x_m,c}$ for some $n>0$.
  Therefore, it reduces to $0$ modulo $(s)$.
\end{proof}

\begin{cor}\label{cor:central_fiber_reduced}
  $\wP_0$ is reduced and geometrically reduced.
\end{cor}

\begin{lem}\label{lem:theta_n}
  For each $n$ only finitely many of the elements
  $\xi_{x+c}\xi_c^{-1}$ in $R_0(c)$ are not zero modulo
  $I^{n+1}=(s^{n+1})$. For $n=0$ the only ones not zero correspond to
  the lattice points $x+c\in\St(\Del_B,c)$, i.e., $x\in\St(\Del_B,0)$.
\end{lem}
\begin{proof}
  Indeed, $\xi_{x+c}\xi_c^{-1}=\zeta_{x,c}s^{B(x,x)/2-B(\alpha,x)}$,
  and $B$ is positive definite. The second part was proved in Lemma
  \ref{lem:paraboloid}.
\end{proof}

\begin{defn}
  We define the Cartier divisor $\wTheta_0$ on $\wP_0$ by the system of
  compatible equations $\{\wtheta/\xi_c \in R_0(c) \}$. Explicitly,
  \begin{displaymath}
    \wtheta/\xi_c = \sum_{x\subset\St(\Del_B,0)}
    a_0(x+c)a^{-1}_0(c) \bar \zeta_{x,c}.
  \end{displaymath}
\end{defn}

\begin{say}
  Clearly, $\wtheta$ defines a $Y$-invariant global section of
  $\wcL_0=\wcL|_{\wP_0}$, so the divisor $\wTheta_0$ is $Y$-invariant.
\end{say}

\begin{lem}
  Once the base change in Assumption \ref{assume:base_change_made} has
  been made, for any further finite base change $S'\to S$, one has
  $\wP'\simeq \wP\underset{S}{\times} S'$. In other words, our
  construction is stable under base change.
\end{lem}
\begin{proof}
  This statement is sufficient to check for the semigroup algebras
  $R(c)$, which is obvious using the basis in
  \ref{thm:first_structure_thm_family}\ref{numi_basis_Rc}.
\end{proof}

\subsection{Case of arbitrary abelian part}
\label{subsec:Case of principal polarization and arbitrary abelian
  part:simp}
\skipline

\begin{saynum}\label{saynum:descr_wP_nontriv_ab_part}
  Most of the statements above transfer to the general case without
  any difficulty. The main difference is that $\wP$ is now fibered
  over $A$ instead of a point, and each $U(\sigma)$, resp.
  $U_0(\sigma)$ is an affine scheme over $A$, resp. $A_0$. One easily
  sees that $\wP$ is isomorphic to the contracted product
  $\wP^r\overset{T}{\times}\wG$ of an $r$-dimensional scheme $\wP^r$
  over $S$ corresponding to the positive definite integral-valued
  bilinear form $B(x,y)$ on the $r$-dimensional lattice $X$ with
  $\wG$. Recall that the contracted product is the quotient of
  $\wP^r\underset{S}{\times}\wG$ by the free action of $T$ with the
  standard action on the first factor and the opposite action on the
  second factor. In the same way, $\wP_0 \simeq
  \wP_0^r\overset{T}{\times}\wG_0$. Moreover, the $\bG_m$-torsor
  $\wdL'$ corresponding to the sheaf $\wcL'\simeq \wcL\otimes
  \cM^{-1}$ is the contracted product of the $\bG_m$-torsor $\wdL^r$
  and $\wG\underset{S}{\times}\bG_{m,S}$, and similarly for the
  central fiber.

  The power series $\wtheta=\sum\xi_x$ defines a $Y$-invariant section
  of $\wcL$.
\end{saynum}

\begin{lem}\label{lem:stratification_wP0}
  $\wP_0$ is a disjoint union of semiabelian varieties $\wG_0(\sigma)$
  which are in 1-to-1 correspondence with the Delaunay cells. One
  has 
  \begin{displaymath}
    1\to \orb(\sigma) \to \wG_0(\sigma) \to A \to 0.
  \end{displaymath}
  The closure $\wV(\sigma)$ of $\wG_0(\sigma)$ is a projective variety 
  $\wV^r(\sigma)\overset{T}{\times}\wG_0$.
  
  For a $0$-dimensional cell $c\in X$ the restriction of $\wtheta$ to
  $\wV(c)\simeq A$ is $\xi_x$ and $(\xi_x)=T_{c^t(x)}(\Theta_A)$. In
  particular, $\wTheta=(\wtheta)$ does not contain any of the strata
  entirely.
\end{lem}
\begin{proof}
  The first part follows at once from Theorem
  \ref{thm:toric_descr_family} by applying the contracted product. The
  second part is obvious because all the other $\xi_x$, $x\ne c$ are
  zero on $\wV(c)$.
\end{proof}

\subsection{Taking the quotient by $Y$}
\label{subsec:Taking the quotient by Y}
\skipline

\begin{lem}
  $\wP$ is a relatively complete model as defined in
  \cite[III.3.1]{FaltingsChai90}.
\end{lem}
\begin{proof}
  We do not even recall the fairly long definition of a relatively
  complete model because most of it formalizes what we already have: a
  scheme $\wP$ locally of finite type over $R$ with an ample sheaf
  $\wcL$, actions of $Y$ and $T$ etc.
  
  There are two additional conditions which we have not described yet.
  The first one is the completeness condition. It is quite tricky but
  it is used in \cite{Mumford_AnalyticDegsAVs,FaltingsChai90} only to
  prove that every irreducible component of $\wP_0$ is proper over
  $k$. We already know this from \ref{thm:toric_descr_family}.
  
  The second condition is that we should have an embedding
  $\wG\into\wP$. It is sufficient to give such an embedding for the
  toric case, since then we simply apply the contracted product.  In
  the toric language, $T$ corresponds to a fan in
  $N_{\bR}=\bX_{\bR}^*\oplus\bR$ consisting of the ray
  $\bR_{\ge0}f_0$. A map from this fan to the fan $\Delta$, which
  sends $f_0$ to $\big(1,-dA(\alpha(\sigma))\big)$ for an arbitrary
  maximal-dimensional Delaunay cell $\sigma$, defines an embedding
  $T\into \wP^r$. We have used the fact that $dA(\alpha(\sigma))$ is
  integral here.
\end{proof}

\begin{rem}
  Note that the embedding $\wG\into\wP$ defines a section of $\wP$
  which has absolutely nothing to do with the zero section $z_{\eta}$
  of $A_{\eta}$ and its closure $z$. The embedding $z\into P$ is
  described by the embedding of fans $(\bZ,\bR_{\ge0}f_0)\into
  (N,\Delta)$, $f_0\mapsto f_0$. From this, we see that
  $z_0\in\orb(\sigma)$, where $\sigma$ is the ``bottom'' face of the
  hyperboloid $Q$. It need not be maximal-dimensional.
\end{rem}

\begin{saynum}
  We can now apply Mumford's construction as described in
  \cite{Mumford_AnalyticDegsAVs,FaltingsChai90}. This consists of
  considering all fattenings
  $(\wP_n,\wcL_n)=(\wP,\wcL)\underset{R}{\times} R/I^{n+1}$, their
  quotients $(P_n,\cL_n)=(\wP_n,\wcL_n)/Y$ and then algebraizing this
  system to a projective scheme $(P,\cL)/S$ such that the generic
  fiber $P_{\eta}$ is abelian and $\cL_{\eta}$ defines a principal
  polarization. By \ref{thm:Faltings_Chai_equivalence},
  $(A_{\eta},\cO(\Theta_{\eta}))\simeq (P_{\eta},\cL_{\eta})$, and
  since the polarization is principal, this isomorphism is uniquely
  defined. Thus, we have obtained the extended family.

  To this construction we will add a theta divisor. By Lemma
  \ref{lem:theta_n} for each $n$ the power series $\wtheta$ defines a
  finite sum in each $R_n(c)$. Hence, we have a compatible system of
  $Y$-invariant sections of $\wcL_n$ that descend to compatible
  sections of $\cL_n$ that algebraize to a section $\theta$ of $\cL$.
\end{saynum}

\begin{say}
  From Lemma \ref{lem:stratification_wP0} we have:
\end{say}

\begin{cor}\label{cor:stratification_P0}
  $P_0$ is a disjoint union of semiabelian varieties $G_0(\bar\sigma)$
  which are in 1-to-1 correspondence with the classes of Delaunay
  cells modulo $Y$-translations.  $\Theta_0=(\theta_0)$ does not
  contain any of the strata entirely.
\end{cor}

\begin{defn}
  In the split case, a {\em stable quasiabelian pair}, SQAP for short,
  over a field $k$ is a pair $(P_0,\Theta_0)$ from our construction.
  In general, a pair of a reduced projective variety and an ample
  Cartier divisor over $k$ is called a stable quasiabelian pair if it
  becomes one after a field extension. We will call $P_0$ itself a
  stable quasiabelian variety, SQAV for short.
\end{defn}

\begin{say}
  Let $\Gamma^2(X^*)$ be the lattice of integral-valued symmetric
  bilinear forms on $X\times X$. For each Delaunay decomposition
  $\Del$ let $K(\Del)$ be the subgroup of $\Gamma^2(X^*)$ generated by
  the positive-definite forms $B$ with $\Del_B=\Del$, and set
  $N(\Del)=\Gamma^2(X^*)/K(\Del)$. 
\end{say}

\begin{thm}
  Each SQAP over an algebraically closed field $k$ is uniquely defined
  by the following data:
  \begin{enumerate}
  \item a semiabelian variety $G_0$, $1\to T_0\to G_0 \to A_0\to 0$, 
  \item a principal polarization $\lambda_{A_0}:A_0\to A_0^t$,
  \item a Delaunay decomposition $\Del$ on $X$, where $X=Y$ is the
    group of characters of $T_0$,
   \item a class of bilinear symmetric sections $\tau_0$ of 
     $(c^t_0\times c_0)^*\cP^{-1}_{A_0}$ on $X\times X$
     modulo the action of the torus $K(\Del)\otimes\bG_{m,k}$.
  \end{enumerate}
\end{thm}
\begin{proof}
  It is clear from the construction that an SQAP depends only on
  $A_0,\Theta_0$ and $A,B,\psi_0,\tau_0$.  Since we do not care about
  the origin and the polarization is principal, giving the pair
  $(A_0,\Theta_0)$ is the same as giving the pair
  $(A_0,\lambda_{A_0})$.  Replacing $\psi$ by $\psi_1$ with the same
  homogeneous part does not change the isomorphism classes of
  subalgebras $\Rp$ and $R$, since the relations between $S_y(\cM_0)$
  remain the same.

  The only information we are getting from $B$ is the Delaunay
  decomposition. Moreover, for a fixed $B$ if we replace the
  uniformizing parameter $s$ by $\mu s$, $\mu\in R\setminus I$, then
  the central fiber will not change, but $b(y,x)$ will change to
  $b(y,x) \mu_0^{B(y,x)}$. Therefore we only have the equivalence class
  by the $K(\Del)\otimes\bG_{m,k}$-action.
\end{proof}

\begin{rem}
  We can further divide the data above by the finite group of
  automorphisms of $(A_0,\lambda_{A_0})$ extended by the finite
  subgroup of $\GL(X)$ preserving $\Del$. The quotient data exactly
  corresponds to the $k$-points of the Voronoi compactification of
  $A_g$. Hence, every $k$-point of
  $\overline{A}_g^{\operatorname{VOR}}$ defines a unique SQAP over
  $k$.
\end{rem}

\begin{rem}
  In \cite{Mumford_KodairaDimSiegel} Mumford considered the
  first-order degenerations of abelian varieties over $\bC$. These are
  exactly our pairs in the case where the toric part of $\wG$ is
  1-dimensional, i.e., $r=1$.
\end{rem}

\begin{say}
  In conclusion, we would like to make the following obvious
  observation. 
\end{say}

\begin{lem}
  The family $P\to S$ is flat. The family of theta divisors $\Theta\to
  S$ is also flat.
\end{lem}
\begin{proof}
  Indeed, $S$ is integral and regular of dimension 1, and $P$ is
  reduced and irreducible so the statement follows, e.g., by
  \cite[III.9.7]{Hartshorne77}. The family of divisors $\Theta\to S$
  is flat because $\Theta\cdot P_t$ is defined at every point $t\in S$
  (\cite[III.9.8.5]{Hartshorne77}).
\end{proof}


\section{Further properties of SQAVs}
\label{sec:Geometric properties of stable quasiabelian varieties}

\begin{say}
  All statements in this section are stable under field extensions.
  Therefore without loss of generality we may assume that $k$ is
  algebraically closed. 
\end{say}

\begin{lem}
  $P_0$ is Gorenstein.
\end{lem}
\begin{proof}
  A Noetherian local ring is Gorenstein if and only if its formal
  completion is such (\cite[18.3]{Matsumura89}). Therefore, we can
  check this property on an \'etale cover $\wP_0$ of $P_0$. Moreover,
  the purely toric case suffices, since $\wP_0$ is a fibration over a
  smooth variety $A_0$ with a fiber $\wP_0^r$, locally trivial in
  \'etale topology.
  
  Recall that the scheme $\wP$, as any torus embedding, is
  Cohen-Macaulay.
  
  $\wP_0\subset\wP$ is the union of divisors $\wV(\sigma)$
  corresponding to the 1-dimensional faces $\Delta(\sigma)$ of the fan
  $\Delta$ (i.e.,  to maximal-dimensional Delaunay cells $\sigma$). The
  following is a basic formula for the dualizing sheaf of a torus
  embedding:
  \begin{displaymath}
    \omega_{\wP}= \cO(-\sum \wV(\sigma)) = \cO(-\wP_0).
  \end{displaymath}
  Since $\wP_0$ is Cartier, $\omega_{\wP}$ is locally free, so $\wP$
  is Gorenstein. The scheme $\wP_0$ is then Gorenstein as a subscheme
  of a Gorenstein scheme that is defined by one regular element.
  
  Alternatively, $\wP_0$ is Gorenstein as the complement of the main
  torus in a toric variety, see \cite[p.126, Ishida's
  criterion]{Oda_ConvexBodies}.
\end{proof}

\begin{lem}
  $\omega_{P_0}\simeq\cO_{P_0}$.
\end{lem}
\begin{proof}
  As above, we have the canonical isomorphism
  \begin{displaymath}
    \omega_{\wP}(\wP_0) \simeq \cO_{\wP},
  \end{displaymath}
  which by adjunction gives $\omega_{\wP_0}\simeq\cO_{\wP_0}$. Since
  both sides are invariant under the action of lattice $Y$, this
  isomorphism descends to $P_0$.
\end{proof}

\begin{thm}
  $$h^i(P_0,\cO)=\binom gi$$.
\end{thm}
\begin{proof}
  We want to exploit the fact that $\wP_0$ is built of ``blocks''
  $\wV(\sigma)$ and that the cohomologies of each block are easily
  computable. Indeed, since the fibers of $\wV(\sigma)\to A_0$ are
  toric varieties $\wV^r(\sigma)$ and $H^i(\wV^r(\sigma),\cO)=0$ for
  $i>0$, we have $R^i\pi_*\cO_{\wV(\sigma)}=0$ for $i>0$ and
  $H^i(\wV(\sigma),\cO)=H^i(A_0,\cO)$. It is well-known that for an
  abelian variety $A_0$ these groups have dimension $\binom ai$.
  
  The following important sequence for a union of torus orbits is
  contained in \cite[p.126]{Oda_ConvexBodies}, where it is called 
  Ishida's complex:
  \begin{displaymath}
    0\to \cO_{\wP^r_0}\to 
    \oplus_{\dim\sigma=g} \cO_{\wV^r(\sigma)} \to \dots
    \oplus_{\dim\sigma=0} \cO_{\wV^r(\sigma)} \to 0.
  \end{displaymath}
  The morphisms in this sequence are the restrictions for all pairs
  $\sigma_1\supset\sigma_2$, taken with $\pm$ depending on a chosen
  orientation of the cells.
  
  By taking the contracted product, we obtain a similar resolution
  for the sheaf $\cO_{\wP_0}$, with $\wV^r(\sigma)$ replaced by
  $\wV(\sigma)$. Finally, by dividing this resolution by the
  $Y$-action, we obtain a resolution of $\cO_{P_0}$. In this
  resolution the morphism $\cO_{\wV(\bar\sigma_1)} \to
  \cO_{\wV(\bar\sigma_2)}$ is a linear combination of several
  restriction maps according to the ways representatives of
  $\bar\sigma_1$ contain representatives of $\sigma_2$.  We can now
  compute $H^i(\cO_{\wP_0})$ by using the hypercohomologies of the
  above complex.

  First, consider the special case where $r=g$ and $a=0$. In this case
  each $H^0(\wV(\sigma),\cO)$ is 1-dimensional and the higher
  cohomologies vanish. Therefore, $H^i(P_0,\cO)$ are the cohomologies
  of the complex
  \begin{displaymath}
    0\to \oplus_{\dim\bar\sigma=g} H^0(\wV(\bar\sigma),\cO) \to \dots
    \oplus_{\dim\bar\sigma=0} H^0(\wV(\bar\sigma),\cO) \to 0.
  \end{displaymath}
  But this complex computes the cellular cohomologies of the cell
  complex $\Del_B/Y$ whose geometric representation is homeomorphic to
  $\bR^r/\bZ^r$. Hence $h^i=\binom ri=\binom gi$.
  
  In general, we obtain the spectral sequence
  $E_1^{pq}=H^p(\bR^r/\bZ^r, H^q(A_0,\cO)) \follows H^{p+q}(P_0,\cO)$
  degenerating in degree 1.  Therefore,
  \begin{displaymath}
    h^i(P_0,\cO)= \sum_{p+q=i} \binom rp \binom aq 
    = \binom {r+a}i = \binom gi.
  \end{displaymath}
\end{proof}

\begin{thm}
  For every $d>0$ and $i>0$ one has $h^0(P_0,\cL_0^d)=d^g$ and
  $h^i(P_0,\cL_0^d)=0$.
\end{thm}
\begin{proof}
  Twist the above resolution of $\cO_{P_0}$ by $\cL$. Once again, the
  cohomologies of each building block are easy to compute. Indeed, for
  a toric variety $\wV^r(\sigma)$ higher cohomologies of an ample
  sheaf vanish. Moreover, since by \ref{thm:toric_descr_family} the
  pair $(\wV^r(\sigma),\wcL)$ is the toric variety with a linearized
  ample sheaf corresponding to the polytope $\sigma\subset X_{\bR}$,
  $H^0(\wV^r(\sigma),\wcL^d)$ is canonically the direct sum of
  1-dimensional eigenspaces, one for each point $z\in\sigma\cap X$,
  cf., e.g., \cite[Ch.2]{Oda_ConvexBodies}.  Taking into account that
  higher cohomologies of ample sheaves on abelian varieties vanish, we
  see that $H^i(\wV(\sigma),\wcL^d)=0$ for $i>0$ and
  $H^0(\wV(\sigma),\wcL^d) \simeq H^0(\wV^r(\sigma),\wcL^d)\otimes
  H^0(A_0,\cM^d)$. It is well-known that for the abelian variety $A_0$
  the latter cohomology space has dimension $d^a$.
  
  The above decomposition into eigenspaces extends to the
  hypercohomologies and we obtain
  \begin{displaymath}
    H^i(P_0,\cL^d) \simeq H^0(A_0,\cM^d) \otimes
    (\oplus_{\bar z\in X/dX} W^i_{\bar z}),
  \end{displaymath}
  where $W^i_{\bar z}$ is the $\bar z$-eigenspace of the $i$-th
  cohomology of the complex
  \begin{displaymath}
    0\to \oplus_{\dim\bar\sigma=g} H^0(\wV^r(\bar\sigma),\wcL^d) \to \dots
    \oplus_{\dim\bar\sigma=0} H^0(\wV^r(\bar\sigma),\wcL^d) \to 0.
  \end{displaymath}
  Fix a representative $z\in X/n$ of $\bar z$.  Let $\sigma_0$ be the
  minimal cell containing $z$. There is a 1-to-1 correspondence
  between the cells $\sigma\ni z$ and the faces of the dual Voronoi
  cell $\widehat\sigma_0$. Since these are exactly the cells for which
  the $z$-eigenspace in $H^0(\wV(\sigma),\wcL^d)$ are nonzero (and
  1-dimensional), we see that $W^i_z$ computes the cellular cohomology
  $H^i(\widehat\sigma_0,k)$. Since as a topological space
  $\widehat\sigma_0$ is contractible, $\dim W^0_z=1$ and $W^i_z=0$ for
  $i>0$. Therefore, $h^i(P_0,\cL^d)=0$ and
  $h^0(P_0,\cL^d)=h^0(A_0,\cM^d)\cdot |X/dX|=d^ad^r=d^g$.
\end{proof}

\begin{saynum}
  As a consequence of this theorem, we can write down explicitly a
  basis of $H^0(P_0,\cL_0^d)$. First, let us do this for the toric
  case. For each $\bar z\in X/dX$ fix a representative $z\in X/d$ and
  choose $x_1,\dots,x_d$ with $x_1+\dots+x_d=dz$. Consider a power
  series
  \begin{eqnarray*}
    && \wtheta'(x_1,\dots,x_d) = \sum_{y\in Y=X} 
    \xi'(x_1+y, \dots, x_d+y), \\
    && \xi'(x_1+y, \dots, x_d+y) :=
    \prod_{i=1}^d \xi_{x_i+y}.
  \end{eqnarray*}
  This power series is obviously invariant under the $Y$-action on
  $K[\vartheta,w^x \,|\, x\in X]$ and, under the assumption
  \ref{assume:base_change_made}, there is a unique nonnegative integer
  $n=n(d_1,\dots,x_d)$ such that $\xi(x_1+y, \dots,
  x_d+y)=s^{-n}\xi'(x_1+y, \dots, x_d+y)$ defines a nonzero section of
  $H^0(\wP_0,\wcL_0^d)$, so that
  $\wtheta(x_1,\dots,x_d)=s^{-n}\wtheta'(x_1,\dots,x_d)$ defines a
  nonzero $Y$-invariant section of $H^0(\wP_0,\wcL_0^d)$ which
  descends to a nonzero section $\theta(x_1,\dots,x_d)$ of
  $H^0(P_0,\cL_0^d)$.  Choosing another representatives $z$ and
  $x_1,\dots,x_d$ changes these sections by multiplicative constants.
  Therefore, in the cases where we do not care about these constants
  we will write simply $\xi_z$, $\wtheta_z$ and $\theta_{\bar z}$.
  
  Geometrically, $\xi_z$ is represented by a point on the surface of
  the multifaceted paraboloid $Q$ of Figure~3 
  lying over $z$, and $\wtheta_z$ by the sum of countably many such
  points lying over all $z+y$, $y\in Y=X$.

  In general, for any $e\in H^0(A_0,\cM^d)$ we repeat the procedure
  taking instead of $\xi_{x_i+y}$ sections
  $\xi_{x_i+y}(e)=S^*_{x_i+y}(e)$. This way, after fixing a basis
  $\{e_1,\dots, e_{d^a}\}$ of $H^0(A_0,\cM^d)$ we obtain a basis
  $\theta_{\bar z}(e_i)$ of $H^0(P_0,\cL_0^d)$.
\end{saynum}

\begin{say}
  The following is an easy application of Theorem
  \ref{thm:structure_of_central_fiber}.\ref{enumi:basis_R0}:
\end{say}

\begin{lem}\label{lem:thetas_explicitly}
  For each $e\in H^0(A_0,\cM^d)$ the following open sets coincide:
  \begin{displaymath}
    \{ \xi_z (e) \ne0 \} =
    \wG_0(\sigma) \cap \pi^{-1} \{ S^*_{dz}(e)\ne0 \},
  \end{displaymath}
  where $\sigma$ is the Delaunay cell containing $z$ in its interior
  $\sigma^0$.
\end{lem}


\begin{thm}
  Let $(P_0,\cL_0)$ be an SQAV of dimension $g$.  Then the sheaf
  $\cO(\cL_0^d)$ is very ample if $d\ge 2g+1$.
\end{thm}
\begin{proof}
  Consider the toric case first, i.e., assume $r=g$, $a=0$. We need to
  prove:
  \begin{numerate}
  \item Sections $\theta_{\bar z}$, $z\in X/dX$ separate the torus
    orbits $\orb(\bar \sigma)$.
  \item They embed every $\orb(\bar \sigma)$.
  \item This embedding is an immersion at every point.
  \end{numerate}
  Let $\sigma$ be a maximal-dimensional cell and pick a point $z\in
  \sigma^0\cap X/d$ in its interior, which exists by
  \ref{lem:comb_of_star}. By \ref{lem:thetas_explicitly} $\xi_z$ is
  nonzero exactly on $\orb(\sigma)$. Therefore, $\wtheta_z$ is nonzero
  exactly on the union of $Y$-translates of this orbit, and
  $\theta_{\bar z}$ is not zero precisely on $\orb(\bar\sigma)$. Thus,
  we have separated the points of this orbit from all the others.
  Continuing by induction down the dimension, we get (i).
  
  The restriction of each $\theta_{\bar z}$ to $\orb(\bar\sigma)$ is a
  sum of degree $d$ monomial corresponding to $\bar z$, provided $\bar
  z\in\bar \sigma$. If $\bar z\in \bar\sigma^0$, there is just one
  such monomial.  Therefore, the condition that suffices for (ii) is
  that the differences of vectors in $\sigma^0\cap X/d$ generate
  $\bR\sigma\cap X/d$ as a group.  This hold by
  \ref{lem:comb_of_star}.\ref{enumi:differences_generate}.
  
  It suffices to prove the immersion condition for the 0-dimensional
  orbit $p$ corresponding to $0\in X$ only. Indeed, it then holds in
  an open neighborhood which intersects every other orbit, and due to
  the torus action everywhere. Moreover, we can work on the \'etale
  cover $\wP_0$.  We have $\wtheta_0(p)\ne 0$, $\wtheta_z(p)=0$ for
  $\bar z\ne0$, and we want to show that $\wtheta_z/\wtheta_0$
  generate $\fm/\fm^2$, where $\fm$ is the maximal ideal of $R_0(0)$.
  On the other hand, $\xi_0(p)\ne 0$ as well, so we can consider
  $\wtheta_z/\xi_0$ instead.
  
  By \ref{thm:structure_of_central_fiber}\ref{enumi:basis_R0},
  $\fm/\fm^2$ is generated by {\em primitive\/} lattice vectors (see
  definition \ref{defn:star_Delaunay_primitive}).  As an element of
  $R_0(0)$, $\wtheta_z/\xi_0$ is the sum with nonzero coefficients of
  monomials $\bar\zeta_{dz'}$ with $z'\in\St(0)$ and $\bar z'=\bar z$.
  Therefore, for (iii) it suffices to have $\Prim/d\subset \St(0)$ and
  $(\Prim-\Prim)\cap dX = \{0\}$. This follows by
  \ref{lem:comb_of_star} again.
  
  Next, assume that the abelian part $A_0$ is nontrivial. To separate
  the orbits $\wG_0(\bar\sigma)$ and the points in the orbit, and to
  see the injectivity, repeat the above arguments with $\xi_z(e_i)$
  such that $e_i$ provide an embedding of $A_0$.
  
  The immersion is again sufficient to prove at the minimal
  dimensional stratum. For every $p\in A_0\subset P_0$ for the tangent
  space we have $\bT_{p,P_0}\simeq \bT_{p,\wP_0^r}\oplus \bT_{p,A_0}$.
  Hence, if $e_0(p)\ne 0$ then to generate $\fm/\fm^2$ it is
  sufficient to take $\theta_{\bar z}(e_z)/\theta_0(e_0)$ and
  $\theta_0(e_i)/\theta_0(e_0)$, where $e_z\big(T_{c^t(dz)}(p)\big)\ne
  0$ and $e_i/e_0$ generate $\bT^*_{p,A_0}$.
\end{proof}  

\begin{lem}\label{lem:comb_of_star}
  \begin{numerate}
  \item $\Prim\subset r \St(0)$.
  \item Assume $d\ge r+2$. Then for each Delaunay cell $\sigma$ the
    differences of vectors in $\sigma^0\cap X/d$ generate
    $\bR\sigma\cap X/d$.
    \label{enumi:differences_generate}
  \item $\big(\St(0)-\St(0)\big)\cap (2+\varepsilon) X= \{0\}$ for any 
    $0<\varepsilon\ll1$.
  \end{numerate}
\end{lem}
\begin{proof}
  Since the restriction of the Delaunay decomposition to $\bR\sigma$
  is again a Delaunay decomposition, we can assume that $\sigma$ is
  maximal-dimensional.
  
  Let $\sigma\subset\St(0)$ be a Delaunay cell and let
  $w\in\Cone(0,\sigma)$ be a primitive lattice vector.  Choose
  arbitrary $r$ linearly independent Delaunay vectors $v_1,\dots,v_r$
  with $w\in \Cone(v_1\dots v_r)$ and write $w=\sum p_iv_i$ for some
  $p_i\in\bQ$. Then obviously $p_i\le1$ and $w/r$ belongs to the
  convex hull of $0,v_1,\dots,v_r$, which is a part of $\St(0)$.  This
  proves (i).
  
  For (ii) note that the vectors $\sum_{i=1}^r v_i$, $(\sum v_i)+v_j$,
  $j=1\dots r$ and $(\sum v_i)-w$, with $w\ne v_i$ primitive all
  belong to $\big((r+2)\sigma\big)^0$.
  
  For (iii), let $\sigma_1\ne\sigma_2$ be two Delaunay cells in
  $\St(0)$. Then for any $y\in X$ the intersection $\sigma_1\cap T_y
  \sigma_2$ is either $\sigma_1$, or a proper face of $\sigma_1$ or
  empty. In the first case $y$ must be a Delaunay vector, since both
  $\sigma_i$ contain $0$. Therefore, $y\notin 2X$. Consequently, for
  any $y\ne0$, $\St(0)\cap T_{2y}\St(0)$ has no interior, and
  $\St(0)\cap T_{(2+\varepsilon)y}\St(0)=\emptyset$.
\end{proof}

\begin{rem}
  The bound above is certainly not optimal. However, it seems that a
  better bound would require going much deeper into the combinatorics
  of Delaunay cells. In the proof above the only properties we used
  were that $\Del_B$ is $X$-periodic and that a Delaunay cell does not
  contain lattice points except its vertices.
\end{rem}


\section{Additions}
\label{sec:Complex-analytic case}

\subsection{Complex-analytic case}
\label{subsec:Complex-analytic case}
\skipline

\begin{say}
  Everything works the same way as in the algebraic case, only easier.
  The ring $\cR$, resp. $K$ is replaced by the stalk of functions
  homomorphic, resp. meromorphic in a neighborhood of $0$. A major
  simplification comes in the construction of the quotient
  $(\wP,\wcL)/Y$. Considering the fattenings $(\wP_n,\wcL_n)/Y$ and
  then algebraizing is unnecessary, since the $Y$-action is properly
  discontinuous in classical topology over a small neighborhood of
  $0$.  Hence, one can take the quotient directly.
  
  The combinatorial description of the family and the central fiber
  and the data for an SQAP remain the same.
\end{say}

\subsection{Higher degree of polarization}
\label{subsec:Higher degree of polarization}
\skipline

\begin{say}
  The formulas in our construction are set up in such a way that we
  can repeat it for any degree of polarization. The outcome, after a
  finite base change, is a normal family with reduced central fiber
  and a relatively ample divisor. However, in this case there are
  several additional choices to make:
  \begin{enumerate}
  \item an embedding $\cM_x\into \cM_{x,\eta}$ for each nonzero
    representative of $X/Y$,
  \item a section $\theta_{A,x}\in H^0(A,\cM_x)$ for each
    representative $x\in X/Y$.
  \end{enumerate}
  According to the description of $H^0(A_{\eta},\cL_{\eta})$ in
  \cite[II.5.1]{FaltingsChai90}, this data is equivalent to providing
  a theta divisor on the generic fiber.  This is why there are
  infinitely many relatively complete models in the case of higher
  polarization.
\end{say}



\section{Examples}
\label{sec:Examples}

\begin{say}
  Below we list all the SQAPs in dimensions 1 and 2, for illustration
  purposes. They can already be found in in
  \cite{Namikawa_NewCompBoth, Namikawa_ToroidalCompSiegel} (over
  $\bC$).
\end{say}

\subsection{Dimension 1}
\label{sec:Dimension 1}
\skipline

\begin{say}
  In this case there is only one Delaunay decomposition of
  $\bZ\subset\bR$, so there is only one principally polarized stable
  quasiabelian pair besides the elliptic curves. The $0$-dimensional
  Delaunay cells correspond to integers $n$, and $1$-dimensional cells
  to intervals $[n,n+1]$.  By \ref{thm:toric_descr_family} the
  corresponding toric varieties are projective lines $(\bP^1,\cO(1))$
  intersecting at points, and the intersections are transversal by
  \ref{thm:structure_of_central_fiber}. The theta divisor restricted
  to each $\bP^1$ has to be a section of $\cO(1)$, i.e., a point.  The
  quotient $P$ is, obviously, a nodal rational curve.
  
\end{say}

\subsection{Dimension 2}
\label{sec:Dimension 2}
\skipline

\begin{say}
  Let us look at the case of the maximal degeneration first. There are
  only two Delaunay decompositions shown on Figures 1,2. In each case
  the irreducible components are the projective toric varieties
  described by the lattice polytopes $\sigma$, see
  \ref{thm:toric_descr_family}.
  
  In the first case we have a net of $(\bP^1\times\bP^1,\cO(1,1))$'s
  intersecting transversally at lines which in turn intersect in fours
  at points. The quotient by $\bZ^2$ will have one irreducible
  component. It is obtained from $\bP^1\times\bP^1$ by gluing two
  pairs of zero and infinity sections.  Modulo the action of
  $K(\Del)\otimes\bG_m$ we have only one parameter $z=b(e_1,e_2)$ and
  the SQAPs with $z$ and $1/z$ are isomorphic.  Therefore, we have a
  family of SQAPs of this type parameterized by $k^*/\bZ_2=k$. The
  theta divisor is the image of a conic on $\bP^1\times\bP^1$. It is
  reduced. It is irreducible unless $z=1$, in which case it is a pair
  of lines.
  
  In the second case \ref{thm:toric_descr_family} we get a net of
  projective planes $(\bP^2,\cO(1))$ meeting at lines which in turn
  meet at points. The quotient will have two irreducible component,
  since there are two non-equivalent maximal-dimensional Delaunay
  cells modulo the lattice. In this case $K(\Del)\otimes\bG_m$ is
  3-dimensional, so a variety $(P_0,\cL_0)$ of this type is unique up
  to isomorphism. The theta divisor restricted to $\bP^2$ has to be
  a section of $\cO(1)$, i.e., a line.  Therefore, the theta divisor on
  $P$ is a union of two rational curves and it is easy to see that
  they intersect at 3 points, one for each $\bP^1$. In other words,
  the theta divisor in this case is a ``dollar curve''.

\end{say}

\begin{rem}
  Note that the lattices in the above two examples are of types
  $A_1\oplus A_1$ and $A_2$ respectively (see f.e.
  \cite{ConwaySloane93}). The number 4, resp. 6, of branches meeting
  at the $0$-dimensional strata has an interesting interpretation in
  this case. For any lattice of the $A,D,E$-type this is what in the
  lattice theory called the ``kissing number'' (think of the billiard
  balls with centers at the lattice elements, each ball is ``kissed''
  by 4, resp. 6, other balls).
\end{rem}

\begin{say}
  There is only one case for a nontrivial abelian part (besides the
  smooth abelian surfaces): when $A_0$ is an elliptic curve. This is
  the simplest case of what Mumford called ``the first order
  degenerations of abelian varieties'' in
  \cite{Mumford_KodairaDimSiegel}.
  
  Before dividing by the lattice $Y=\bZ$ we have a locally free
  fibration over an elliptic curve with a fiber which corresponds to
  the case of maximal degeneration of dimension 1, i.e,. a chain of
  projective lines. The group $Y$ acts on this scheme $\wP_0$ by
  cycling through the chain and at the same time shifting ``sideways''
  with respect to the elliptic curve. The theta divisor $\wTheta_0$ on
  $\wP_0$ is invariant under this shift, so it descends to a divisor
  on $P_0$.
\end{say}


\ifx\undefined\bysame
\newcommand{\bysame}{\leavevmode\hbox to3em{\hrulefill}\,}
\fi

\end{document}